\documentclass[11pt]{amsart}
\usepackage{mymathmacros}
\usepackage{formatmacros}
\usepackage{analysismacros}
\usepackage{alphamacros}
\usepackage{linalgmacros}
\usepackage{refmacros}
\usepackage{multido}

\usepackage[monochrome]{pstricks}
\usepackage{pst-plot, pst-math, pst-func}
\usepackage{pstricks-add}

\begin{document}

\title[Mendes Conjecture for time-one maps]{Foliations and Conjugacy, II: \\The Mendes Conjecture \\for time-one maps of flows}
\author{Jorge Groisman}
	\address{Instituto de Matem\'atica y Estad\'{\i}stica Prof. Ing. Rafael Laguardia, 
	Facultad de Ingenier\'{\i}a Julio Herrera y Reissig 565 11300, MONTEVIDEO, Uruguay}
	\email{jorgeg@fing.edu.uy}
\author{Zbigniew Nitecki}
	\address{Department of Mathematics, Tufts University, Medford, MA 02155}
	\email{zbigniew.nitecki@tufts.edu}
\thanks{The second author thanks IMERL for its hospitality and support during a visit in August 2018}

\date{}                                           


\begin{abstract}A diffeomorphism \selfmap{f}{\Realstwo} in the plane is \emph{Anosov} if it has a hyperbolic splitting 
at every point of the plane. 
The two known topological conjugacy classes of such diffeomorphisms are linear hyperbolic automorphisms
and translations (the existence of Anosov structures for plane translations was originally shown by W. White).
P. Mendes conjectured that these are the only topological conjugacy classes for Anosov diffeomorphisms in the plane.
We prove that this claim holds when the Anosov diffeomorphism is the time-one map of a flow, via a theorem about 
foliations invariant under a time one map.
\end{abstract}

\newcommand{\clints}[3]{\ensuremath{\clint{#1}{#2}_{#3}}}
\newcommand{\opints}[3]{\ensuremath{\opint{#1}{#2}_{#3}}}
\newcommand{\ropints}[3]{\ensuremath{\ropint{#1}{#2}_{#3}}}
\newcommand{\lopints}[3]{\ensuremath{\lopint{#1}{#2}_{#3}}}

\newcommand{\uopint}{\opint{0}{1}}

\newcommand{\cF}{\ensuremath{\mathcal{F}}}
\newcommand{\cFs}[1]{\ensuremath{\cF_{#1}}}
\newcommand{\Fclint}[2]{\ensuremath{\clints{#1}{#2}{\cF}}}
\newcommand{\Fopint}[2]{\ensuremath{\opints{#1}{#2}{\cF}}}
\newcommand{\Fropint}[2]{\ensuremath{\ropints{#1}{#2}{\cF}}}
\newcommand{\Flopint}[2]{\ensuremath{\lopints{#1}{#2}{\cF}}}
\newcommand{\Farc}{\cF-arc}
\newcommand{\Fbox}{\cF-box}
\newcommand{\FR}{\ensuremath{{\mathcal{R}_{\cF}}}}

\newcommand{\cFst}{\ensuremath{\cF^{s}}}
\newcommand{\Fsclint}[2]{\ensuremath{\clints{#1}{#2}{\cFst}}}
\newcommand{\cFut}{\ensuremath{\cF^{u}}}
\newcommand{\Fuclint}[2]{\ensuremath{\clints{#1}{#2}{\cFut}}}
\newcommand{\cR}{\ensuremath{\mathcal{R}}}

\newcommand{\Esa}{\ensuremath{E^{s}_{1}}}
\newcommand{\Esb}{\ensuremath{E^{s}_{2}}}
\newcommand{\Eua}{\ensuremath{E^{u}_{1}}}
\newcommand{\Eub}{\ensuremath{E^{u}_{2}}}

\newcommand{\cG}{\ensuremath{\mathcal{G}}}
\newcommand{\cGs}[1]{\ensuremath{\cG_{#1}}}
\newcommand{\Gclint}[2]{\ensuremath{\clints{#1}{#2}{\cG}}}
\newcommand{\Gopint}[2]{\ensuremath{\opints{#1}{#2}{\cG}}}
\newcommand{\Gropint}[2]{\ensuremath{\ropints{#1}{#2}{\cG}}}
\newcommand{\Glopint}[2]{\ensuremath{\lopints{#1}{#2}{\cG}}}
\newcommand{\Glen}{\flow-length}
\newcommand{\Garc}{\cG-arc}
\newcommand{\Gbox}{\cG-box}
\newcommand{\Gline}{\cG-leaf}
\newcommand{\cGso}{\cGs{\xs{0}}}
\newcommand{\Gnbhd}{\cG-neighborhood}
\newcommand{\Gless}{\ensuremath{\sqsubset}}

\renewcommand{\H}{\ensuremath{\mathcal{H}}}
\newcommand{\sideof}[2]{\ensuremath{\H_{#2}^{#1}}}
\newcommand{\sidepos}[1]{\sideof{+}{#1}}
\newcommand{\sideneg}[1]{\sideof{-}{#1}}

\newcommand{\flow}{\vphi}
\newcommand{\flowto}[1]{\ensuremath{\flow^{#1}}}
\newcommand{\flowtoof}[2]{\ensuremath{\flowto{#1}(#2)}}
\newcommand{\flowt}{\flowto{t}}
\newcommand{\flowclint}[2]{\ensuremath{\clints{#1}{#2}{\flow}}}
\newcommand{\flowopint}[2]{\ensuremath{\opints{#1}{#2}{\flow}}}
\newcommand{\flowropint}[2]{\ensuremath{\ropints{#1}{#2}{\flow}}}
\newcommand{\flowlopint}[2]{\ensuremath{\lopints{#1}{#2}{\flow}}}

\newcommand{\J}{\ensuremath{J}}
\newcommand{\prolp}[2]{\ensuremath{\J_{#1}^{+}(#2)}}
\newcommand{\prolpfl}[1]{\prolp{\flowt}{#1}}
\newcommand{\prolpf}[1]{\prolp{f}{#1}}
\newcommand{\proln}[2]{\ensuremath{\J_{#1}^{-}(#2)}}
\newcommand{\prolnfl}[1]{\proln{\flowt}{#1}}
\newcommand{\prolnf}[1]{\proln{f}{#1}}

\newcommand{\Orb}{\ensuremath{\mathcal{O}}}
\newcommand{\orb}[2]{\ensuremath{\Orb_{#1}(#2)}}
\newcommand{\orbfl}[1]{\orb{\flow}{#1}}
\newcommand{\orbf}[1]{\orb{f}{#1}}

\newcommand{\Gamp}{\ensuremath{\Gamma_{+}}}
\newcommand{\Gamn}{\ensuremath{\Gamma_{-}}}
\newcommand{\U}{\ensuremath{\mathcal{U}}}
\newcommand{\bdU}{\ensuremath{\partial\U}}

\newcommand{\Dm}{\ensuremath{D_{-}}}

\newcommand{\tsi}{\ts{i}}
\newcommand{\ksi}{\ks{i}}
\newcommand{\xpps}[1]{\ensuremath{x\ppr_{#1}}}
\newcommand{\ypps}[1]{\ensuremath{y\ppr_{#1}}}
\newcommand{\xppsi}{\ensuremath{\xpps{i}}}
\newcommand{\yppsi}{\ensuremath{\ypps{i}}}
\newcommand{\Tp}{\ensuremath{T\pr}}
\renewcommand{\Tps}[1]{\ensuremath{\Tp_{#1}}}
\newcommand{\Tpp}{\ensuremath{T\ppr}}

\newcommand{\cL}{\ensuremath{\mathcal{L}}}
\newcommand{\cLs}[1]{\ensuremath{\cL_{#1}}}
\newcommand{\Lclint}[2]{\ensuremath{\clints{#1}{#2}{\cL}}}
\newcommand{\Lopint}[2]{\ensuremath{\opints{#1}{#2}{\cL}}}
\newcommand{\Lropint}[2]{\ensuremath{\ropints{#1}{#2}{\cL}}}
\newcommand{\Llopint}[2]{\ensuremath{\lopints{#1}{#2}{\cL}}}
\newcommand{\Larc}{\cL-arc}
\newcommand{\Lbox}{\cL-box}

\newcommand{\Phitoof}[2]{\ensuremath{\Phi^{#1}(#2)}}

\newcommand{\V}{\ensuremath{V}}
\newcommand{\VG}{\ensuremath{\V_{\cG}}}
\newcommand{\VL}{\ensuremath{\V_{\A}}}
\renewcommand{\half}{\ensuremath{\frac{1}{2}}}

\newcommand{\A}{\ensuremath{L_{A}}}

\newcommand{\W}{\ensuremath{\mathcal{W}}}
\newcommand{\Wsof}[1]{\ensuremath{\of{\W^{s}}{#1}}}
\newcommand{\Wuof}[1]{\ensuremath{\of{\W^{u}}{#1}}}

\newcommand{\cJ}{\ensuremath{\mathcal{J}}}
\newcommand{\cJsof}[1]{\ensuremath{\cJ^{s}_{#1}}}
\newcommand{\cJuof}[1]{\ensuremath{\cJ^{u}_{#1}}}
\newcommand{\JAsof}[1]{\ensuremath{{\cJ_{\A}^{s}}_{#1}}}
\newcommand{\JAuof}[1]{\ensuremath{{\cJ_{A}^{u}}_{#1}}}

\newcommand{\Vso}{\Vs{0}}
\newcommand{\tV}{\ensuremath{\tilde{V}}}
\newcommand{\tVso}{\ensuremath{\tV_{0}}}

\renewcommand{\Qs}[1]{\ensuremath{Q_{#1}}}
\newcommand{\Qsi}{\Qs{i}}
\newcommand{\pQsi}{\ensuremath{Q\pr_{i}}}
\newcommand{\tQsi}{\ensuremath{\tilde{Q}_{i}}}
\newcommand{\tg}{\ensuremath{\tilde{g}}}

\renewcommand{\gams}[1]{\ensuremath{\gamma_{#1}}}

\newcommand{\qsi}{\qs{i}}

\newcommand{\metric}{\ensuremath{\mu}}

\newtheorem*{Mendescon}{Mendes' Conjecture}
\newtheorem*{ThmA}{Theorem A}
\newtheorem*{ThmB}{Theorem B}
\newtheorem*{ThmC}{Theorem C}

\maketitle
\section{Introduction}\label{sec:intro}
	A diffeomorphism \selfmap{f}{M} of a compact manifold $M$ is called \deffont{Anosov} 
	if it has a global hyperbolic splitting of the tangent bundle.  
	Such diffeomorphisms have been studied extensively in the past fifty years.
	The existence of a splitting implies the existence of two foliations, into stable \resp{unstable} manifolds, preserved by the
	diffeomorphism, such that the map shrinks distances along the stable leaves, while its inverse does so for the unstable ones.
	Anosov diffeomorphisms of compact manifolds have strong recurrence properties.
	
	The existence of an Anosov structure when $M$ is compact is independent of the Riemann metric used to define it, 
	and the foliations are invariants of topological conjugacy.
	By contrast, an Anosov structure on a non-compact manifold is highly dependent on the Riemann metric, 
	and the recurrence properties observed in the compact case do not hold in general.
	This is strikingly illustrated by Warren White's example \cite{White} of a complete Riemann metric on the plane \Realstwo{}
	for which the horizontal translation is Anosov.  Furthermore, as we showed in an earlier paper \cite{GroismanNi1}, the
	stable and unstable foliations are not invariants of topological conjugacy among Anosov diffeomorphisms.
	
	Prompted by White's example, Pedro Mendes \cite{Mendes} formulated the following
\begin{definition}\label{dfn:Anosov}
	An \deffont{Anosov structure} on \Realstwo{} for a  diffeomorphism \selfmap{f}{\Realstwo} consists of  
	a complete Riemannian metric 
	on \Realstwo{} and 
	\begin{description}
		\item[Stable and Unstable Foliations] two continuous foliations \cFst{} and \cFut{} with \Cr{1} leaves
		varying continuously in the \Cr{1} topology and
		respected by $f$: the image of a leaf of \cFst{} \resp{\cFut} is again a leaf of \cFst{} \resp{\cFut};
		\item[Hyperbolicity] 
		there exist constants $C$ and  $\lam>1$ such that for any positive integer $n$
		and any vector \vv{} tangent to a
		leaf of \cFut{}.
		\begin{equation*}
			\norm{Df^{n}(\vv)}\geq C\lam^{n}\norm{\vv}
		\end{equation*}
		while for any vector \vv{} tangent to a leaf of \cFst{}
		\begin{equation*}
			\norm{Df^{n}(\vv)}\leq C\lam^{-n}\norm{\vv}
		\end{equation*}
		where \norm{\vv} denotes the length of a vector using the metric \metric.
	\end{description}
\end{definition} 
We shall use the adjectives \emph{Anosov}, \emph{stable} and \emph{unstable} in the natural way:  
 a diffeomorphism is \deffont{Anosov} if it has an  Anosov structure; the leaf of \cFst{} \resp{\cFut} through a point is
 its \deffont{stable} \resp{\deffont{unstable}} \deffont{leaf}.  
 
 He proved several general properties of Anosov diffeomorphisms of the plane, and asked if the two known examples
 represent all possible topological conjugacy classes among them:  
 
\begin{Mendescon}
	If an orientation-preserving diffeomorphism \selfmap{f}{\Realstwo} has an Anosov structure, then $f $ is topologically conjugate to either 
	\begin{itemize}
		\item the translation 
		\begin{equation*}
			T(x,y)= (x+1,y)
		\end{equation*}
		 or 
		 \item the hyperbolic linear automorphism \selfmap{\A}{\Realstwo} defined by
		 \begin{equation*}
		 	\of{\A}{\vx}=A\vx
		\end{equation*}
		where 
		\begin{equation*}
			A=\left[\begin{array}{cc}2 & 0 \\0 & \frac{1}{2}\end{array}\right].
		\end{equation*}
	\end{itemize}
\end{Mendescon}

In a first step toward establishing this conjecture, Mendes proved  
\begin{theorem}[Mendes, \cite{Mendes}]\label{thm:MendesThm}
	If \selfmap{f}{\Realstwo} is a diffeomorphism of the plane with an Anosov structure, then
	\begin{enumerate}
		\item $f$ has at most one nonwandering point (which then must be a hyperbolic fixedpoint);
		\item any point with nonempty $\alpha$-\resp{$\omega$-}limit set 
		has empty forward \resp{backward} prolongational limit set under $f$.
	\end{enumerate}
\end{theorem}  

In this paper, we establish the truth of Mendes' conjecture under an additional assumption:
\begin{ThmA}
	If \flowt{} is a \Cr{1} flow on \Realstwo{} and $f=\flowto{1}$ is its time-one map, then the existence of an Anosov structure
	for $f$ implies the conclusion of Mendes' Conjecture.
\end{ThmA}

Our proof divides into the two cases given by the first conclusion in \refer{thm}{MendesThm}:
\begin{description}
	\item[Case 1] $f$ has empty nonwandering set (\ie{} $f$ is is a ``Brouwer translation'');
	\item[Case 2] $f$ has a unique nonwandering point.
\end{description}

In the first case, the assumption that $f=\flowto{1}$ is fixedpoint-free implies that the flow \flowt{} has no fixedpoints.
Thus the flowlines of \flowt{} form a foliation \cG{} of \Realstwo.  
A foliation \cG{} of \Realstwo{} of  is \deffont{trivial} if there is a homeomorphism \selfmap{H}{\Realstwo} taking leaves
of \cG{} to horizontal lines. 
Triviality of the orbit foliation of a fixedpoint-free flow is equivalent to topological conjugacy of its time-one map 
with a translation (\refer{prop}{trivial}).

In \refer{subsec}{action}, we establish the following theorem about 
foliations preserved by the time-one map of a nontrivial flow:%
\footnote{
	This does not assume an Anosov structure for $f$.
}

\begin{ThmB}
	Suppose \flowt{} is a fixedpoint-free \Cr{1} flow in the plane with flow line foliation \cG.  
	Let \selfmap{f=\flowto{1}}{\Realstwo} be the time-one map
	of \flowt, and suppose \cF{} is a \Cr{1} foliation preserved by $f$.
	
	If \cG{} is nontrivial, then some leaf of \cF{} is invariant under $f$.
\end{ThmB} 

In \refer{subsec}{trivReeb}, we give a proof of Theorem C, a characterization of nontrivial foliations in terms of the existence of 
nontrivial prolongation relations between leaves (``Reeb components''), which forms the basis of our proof of Theorem B.

Applying this to the stable foliation in case $f$ is Anosov, we see that when the time-one map of a flow on \Realstwo{} is an Anosov
Brouwer translation, the orbit foliation must be trivial, and hence the map must be topologically conjugate to a translation
(\refer{cor}{Mendes1}).

In the second case, the unique nonwandering point of $f$ must be a fixedpoint of the flow;%
\footnote{
The only way a fixedpont of a time-one map is not a fixedpoint of the flow is if it lies on a period-one closed orbit
of the flow.  By the Poincar\'e-Bendixson Theorem, this would force a fixedpoint of the flow elsewhere in the plane.
}
 the presence of an Anosov structure means
that it is a hyperbolic saddle point, and the stable \resp{unstable} leaf through this point 
consists of the two incoming \resp{outgoing} separatrices together with the fixedpoint itself.
The second conclusion in \refer{thm}{MendesThm} implies that these separatrices escape to infiinity, and hence
separate the plane into four quadrants.  
Another application of Theorem A to the restriction of $f$ to any one of these quadrants and to its stable foliation 
shows that the restriction of the foliation to each (open) quadrant is trivial.
In \refer{subsec}{fixedpoint} we use a standard ``fundamental domain'' argument to construct a topological conjugacy 
between the restrictions of $f$ and \A{} to  invariant neighborhoods of the fixedpoint, and then use the triviality of the 
flow line foliation \cG{} in each quadrant to extend this conjugacy to the whole plane.

A subtle point here is that examples in \cite{GroismanNi1} show that in general we cannot hope to preserve the stable and unstable foliations under
this conjugacy.  
Although the conjugacy we construct on the invariant neighborhood of the fixedpoint \emph{does} preserve the restriction of
these foliations to the neighborhood, the extension to the rest of the plane need not do so.

\section{Foliations Invariant under a time one map}\label{sec:FITOM}

In this section, we prove Theorem B.
%

\subsection{Preliminaries on plane foliations}\label{subsec:PlanFolPrelims}

A nonvanishing \Cr{1} vectorfield in the plane generates a fixedpont-free flow \flowt{} 
whose (directed) orbits form an oriented foliation of \cG{} of \Realstwo; 
conversely every \Cr{1} foliation \cG{} of \Realstwo{} can be oriented, and viewed as the set of orbits of some 
\Cr{1}
flow on \Realstwo.
We adopt interval notation for arcs in \Realstwo: for example, a closed arc with endpoints $x,y\in\Realstwo$
will be denoted \clint{x}{y}.
Given a point $x\in\Realstwo$, the leaf of \cG{} through $x$ is denoted $\boldsymbol{\cGs{x}}$; if $\xp\in\cGs{x}$, 
the closed arc of \cG{} joining $x$ and \xp{} is denoted $\boldsymbol{\Gclint{x}{\xp}}$  
When \cG{} is the foliation by orbits of the flow \flowt,
we can write $\boldsymbol{\flowclint{x}{x+t}}$ for \Gclint{x}{\flowt(x)}; in this case we can assign to any \Garc{} \Gclint{a}{b}
the \deffont{\Glen} \abs{t-s}, where $a=\flowto{s}(x)$ and $b=\flowto{t}(x)$.

We can extend these ideas and notations to open \resp{half-open} \cG-arcs \Gopint{x}{\xp}
 \resp{\Gropint{x}{\xp}, \Glopint{x}{\xp}}.


An arc is (topoogically \deffont{transverse} to the foliation \cG{} if it crosses each leaf of \cG{} at most once.
A closed \resp{open} \deffont{$\boldsymbol{\cG}$-box} is a topological disc $D\subset\Realstwo$ homeomorphic to the
unit square $\uclint\times\uclint$\resp{$\uopint\times\uopint$} such that the horizontal \resp{vertical} arcs in the square correspond to \Garc{}s 
\resp{transversals to \cG}.

\begin{remark}\label{rmk:Gbox}
	Given any \Garc{} $\gamma=\Gclint{\xs{0}}{\xs{1}}$ and transversals to \cG{} \Ts{i} at \xs{i}, $i=0,1$, 
	we can find a \Gbox{} containing \gam{} whose vertical edges are subtransversals of \Ts{i} at \xs{i}.
	We call this a $\boldsymbol{\Gnbhd}$ of \gam{}.
	This last idea can be extended to a full leaf: given a transversal $T$ at $x$, the union of all leaves through $T$
	is naturally homeomorphic to the cartesian product $T\times\cGs{x}$ and contains \cGs{x} in its interior;
	we call this a \deffont{tubular neighborhood} of the leaf \cGs{x}. 
	For any compact arc \Gclint{\xs{0}}{\xs{1}} contained in \cGs{x}, 
	any pair of transversals \Ts{i} at \xs{i} $i=0,1$ cuts off a \Gnbhd{} of the arc;
	note in particular that there are no recurrent leaves: the intersection of the full leaf with this neighborhood 
	consists of the arc alone.  The arc \gam{} separates the \Gnbhd{} into two \textbf{one-sided $\boldsymbol{\Gnbhd}$s} 
	of \gam.
		
\end{remark}

When \cG{} is generated by the flow \flow, we can measure the ``\flow-size'' of a
\Gbox: its \deffont{height} is the the maximum of the lengths of its vertical sides; its \deffont{outer \Glen} 
\resp{\deffont{inner \Glen}} is the maximum \resp{minimum} among the \Glen{s} of its horizontal arcs.

\subsubsection{Trivial foliations}\label{subsec:trivReeb}

The orbits of the translation flow $\flowtoof{t}{(x,y)}=(x+t,y)$ are horizontal lines.  
A plane foliation \cG{} is \deffont{trivial} if there is a homeomorphism of the plane to itself 
taking the leaves of \cG{} to horizontal lines;  
this is equivalent to the existence of a \deffont{global cross section} 
to the foliation--an arc that meets each leaf of \cG{} exactly once.

The obstruction to triviality can be described in the language of prolongations.

\begin{definition}\label{dfn:prolfol}
	Two points $x,y\in\Realstwo$ on distinct leaves of a foliation are \deffont{prolongationally related} in \cG{} 
	if there exist points
	$\xsi\to x$ and $\ysi\to y$ such that for each $i=1,2,...$, \xsi{} and \ysi{} lie on the same leaf of \cG.
\end{definition}
It is easy to see that in such a situation every point of \cGs{x} is prolongationally related to every point of \cGs{y}.

This is essentially related to the \emph{prolongational limit sets} of dynamical systems: 
\begin{definition}\label{dfn:proldyn}
	Given a flow \flowt, the point $y\in \Realstwo$ is in the (first) \deffont{forward prolongation} of $x\in\Realstwo$
	under \flowt, $y\in \boldsymbol{\prolpfl{x}}$, if there exist points $\xsi\to x$ and times $\tsi\to+\infty$ such that 
	$\ysi\eqdef\flowtoof{\tsi}{\xsi}\to y$.
	
	The \deffont{backward prolongation} $\boldsymbol{\prolnfl{x}}$ is defined as above, with $+\infty$ replaced by $-\infty$.
\end{definition}
Clearly, $y\in\prolpfl{x}$ if and only if $x\in\prolnfl{y}$, and each set is invariant under the flow.  (However, in general neither
set consists of a single \flow-orbit.)  Also, when \cG{} is generated by a flow, two points are prolongationally related under \cG{}
if and only if each is in the forward or backward prolongation of the other under \flowt.

Suppose two points $x,y\in\Realstwo$ are prolongationally related in \cG.  Each leaf \cGs{x} \resp{\cGs{y}} separates \Realstwo,
so the complement of their union has three components;  we denote by \U{} the component whose boundary consists
of \emph{both} \cG-lines.

\begin{definition}\label{dfn:Reeb}
	A \deffont{Reeb component} of a foliation of \Realstwo{} is an open topological disc \U{} whose frontier in \Realstwo{}
	consists of two \cG-lines, $\Gamn=\cGs{x}$ and $\Gamp=\cGs{y}$, which are prolongationally related.
\end{definition}
We note that the choice of which edge is labelled \Gamp{} is dictated by an orientation of the foliation;  when this comes from
a flow we say $\Gamn=\cGs{x}$ when $y\in\prolpfl{x}$.


			\begin{figure}[htbp]
			\begin{center}
				\begin{pspicture}(-4.5,-2.5)(2.5,2.5)
					\pspolygon[linestyle=none, fillstyle=solid, fillcolor=gray!25](2.8,-1.8)(-4.5,-1.8)(-4.5,1.8)(2.7,1.8)
				
					\rput*(-3,0){\U}
					\psline[linewidth=1.8pt, arrowscale=1.5]{->}(-4.8,1.8)(2.8,1.8)
						\uput[r](2.8,2){$\Gamma_{-}$}
					\psline[linewidth=1.8pt, arrowscale=1.5]{->}(2.8,-1.8)(-4.8,-1.8)
						\uput[r](2.8,-2){$\Gamma_{+}$}
					\parametricplot[plotstyle=curve, plotpoints=50]{-0.30}{0.30}{0.8 t Pi mul COS div neg -3.0 add t 4 mul}
					\parametricplot[plotstyle=curve, plotpoints=50]{-0.35}{0.35}{0.8 t Pi mul COS div neg -2.5 add t 4 mul}
					\parametricplot[plotstyle=curve, plotpoints=50]{-0.38}{0.38}{0.8 t Pi mul COS div neg -2.0 add t 4 mul}
					\parametricplot[plotstyle=curve, plotpoints=50]{-0.40}{0.40}{0.8 t Pi mul COS div neg -1.5 add t 4 mul}
					\parametricplot[plotstyle=curve, plotpoints=50]{-0.43}{0.43}{0.8 t Pi mul COS div neg -1.0 add t 4 mul}
					\parametricplot[plotstyle=curve, plotpoints=50]{-0.42}{0.42}{0.8 t Pi mul COS div neg -0.5 add t 4 mul}
					\parametricplot[plotstyle=curve, plotpoints=50]{-0.43}{0.43}{0.8 t Pi mul COS div neg t 4 mul}
					\parametricplot[plotstyle=curve, plotpoints=50]{-0.435}{0.435}{0.8 t Pi mul COS div neg 0.5 add t 4 mul}
					\parametricplot[plotstyle=curve, plotpoints=50]{-0.443}{0.443}{0.8 t Pi mul COS div neg 1.0 add t 4 mul}
					\parametricplot[plotstyle=curve, plotpoints=50]{-0.45}{0.45}{0.8 t Pi mul COS div neg 1.5 add t 4 mul}
					\parametricplot[plotstyle=curve, plotpoints=50]{-0.455}{0.455}{0.8 t Pi mul COS div neg 2.0 add t 4 mul}
					\parametricplot[plotstyle=curve, plotpoints=50]{-0.46}{0.46}{0.8 t Pi mul COS div neg 2.5 add t 4 mul}
					\parametricplot[plotstyle=curve, plotpoints=50]{-0.46}{0.46}{0.8 t Pi mul COS div neg 3.0 add t 4 mul}
					\parametricplot[plotstyle=curve, plotpoints=50]{-0.46}{0.46}{0.8 t Pi mul COS div neg 3.5 add t 4 mul}
				\end{pspicture}
			\caption{A Reeb component}
			\label{fig:Reeb}
			\end{center}
			\end{figure}
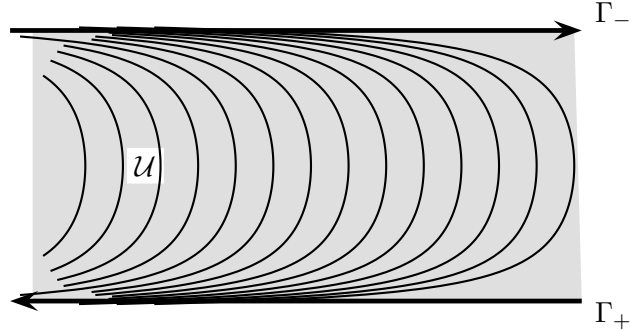

\begin{ThmC}
	A foliation of the plane is trivial if and only if it has no Reeb components.
\end{ThmC}
The basic idea is implicit in the work of Whitney \cite{Whitney}, Kaplan \cite{Kaplan1}, 
and Haefliger, Reeb, and Godbillon \cite{Haefliger-Reeb},\cite{Godbillon-Reeb},\cite{Godbillon}.  
However, none of these works gives an explicit statement of this equivalence, 
even without the language of prolongational limits.  

A set $\U\subset \Realstwo$ is \deffont{saturated} by the foliation \cG{}  if it is a union of leaves;
when \cG{} is an orbit foliation, this is equivalent to invariance of \U{} under the flow.
We say a connected, saturated set \U{} has the \deffont{separation} property if, given any three leaves in \U{}
one of them separates the other two.

When \U{} has the separation property, the leaves in \U{} can be linearly ordered as follows:
\begin{enumerate}
	\item Pick a ``base leaf'' \cGso{} and designate one of the complementary (topological) half-planes 
	as the ``positive'' side \sidepos{\xso} of \cGso{} and the other as its ``negative'' side \sideneg{\xso}.
	\item If the leaf \cGs{x} is in the positive half-plane \sidepos{\xso}, then its ``positive'' side \sidepos{x} 
	is the one disjoint from \cGso{};
	if \cGs{x} is in the ``negative'' half-plane \sideneg{\xso}, then its ``negative'' side \sideneg{x} is the one disjoint from \cGso.
	\item With these designations, write $\boldsymbol{\cGs{x}\Gless\cGs{y}}$ if \cGs{x} is on the negative side \sideneg{y}
	of \cGs{y}
	(equivalently, if \cGs{y} is on the positive side \sidepos{x} of \cGs{x}).
\end{enumerate}
We can, by abuse of notation, write $x\leq y$ for any two points of \U, with the understanding that ``$x<y$''
means $\cGs{x}\Gless\cGs{y}$ but ``$x=y$ '' means only that they lie on the same leaf.

If \U{} has a global cross section, then this ordering of leaves corresponds to the order of their intersections 
with the cross section.  This gives one direction of the following observation, implicit in \cite{Kaplan1}:
\begin{prop}\label{prop:sep}
	An open, connected, saturated set has a global cross section if and only if it has the separation property.
\end{prop}

	In view of the preceding observation, to prove \refer{prop}{sep} we need only show that 
	every connected, open saturated set with separation has a global cross section.
	
	We will do this by piecing together local cross sections, using the following observation.
	
\begin{remark}\label{rmk:transversals}
	\begin{enumerate}
		\item If \clint{x}{y} and \clint{y}{z} are transversal arcs separated by \cGs{y}, then their union 
		\clint{x}{z} is transversal to \cG{}.
		\item If \clint{x}{\ys{1}} and \clint{\ys{2}}{z} are transverse arcs such that \ys{1} and \ys{2} lie on the same
		leaf \cGs{y}, which separates these arcs, then there is a transverse arc \clint{x}{z} which agrees with their
		union except for an arbitrarily chosen \Gnbhd{} of \Gclint{\ys{1}}{\ys{2}}.
		In fact, this arc can be chosen to cross \cGs{y} at any desired point on \Gclint{\ys{1}}{\ys{2}},
		including the possibility that \clint{x}{z} contains one of \clint{x}{\ys{1}} or \clint{\ys{2}}{z}.
	\end{enumerate}
\end{remark}

	Fix $\U\subset\Realstwo$ an open, connected, saturated set with separation.
	First, we establish
	\begin{lemma}\label{lem:connect}
		Every pair of points $x,y\in\U$ lying on distinct \cG-leaves can be joined by a transverse arc
		$\clint{x}{y}\subset\U$.
	\end{lemma}
	\begin{proofof}{\refer{lem}{connect}}
		Given two leaves $\cGs{x}\Gless\cGs{y}$ in \U{} which can be joined by a transverse arc \clint{x}{y}, 
		the set of all leaves crossing \opint{x}{y}--equivalently the set of leaves \cGs{z} satisfying
		$\cGs{x}\Gless\cGs{z}\Gless\cGs{y}$--is the \deffont{open strip} with edges \cGs{x} and \cGs{y}.
		These strips cover \U{}, and clearly any pair of points on distinct leaves but in the same strip 
		can be joined by a transversal arc contained in the strip.  
		
		A compact arc \clint{x}{x\pr} can be partitioned into finitely many successive subarcs \clint{\xs{i}}{\xs{i+1}}
		in such a way that the endpoints of each subarc belong to a common strip; of course there
		may be more than one such pair (not necessarily for adjacent subarcs) belonging to the same strip.
		
		Now given a curve \gam{} intersecting every leaf in \U{} (not necessarily transversally) we can partition it into
		compact arcs and hence by the observation above, we can find a bisequence of points 
		\xs{i}, $i\in\Integers$ 
		partitioning \gam{} such that for each $i\in\Integers$. \xs{i} and \xs{i+1} can be joined by an arc \clint{\xs{i}}{\xs{i+1}}
		transverse to \cG.  These arcs might not form a global cross section, but we can modify the bisequence
		to make it monotone with respect to the \Gless{} relation: fix a ``base'' partition point \xs{\is{0}} and 
		then define \is{k} recursively for $k>0$ by setting $\is{k}=\min\setbld{j>\is{k-1}}{\cGs{\xs{\is{k-1}}}\Gless\cGs{\xs{j}}}$;
		similarly for $k<0$ set $\is{k}=\max\setbld{j<\is{k+1}}{\cGs{\xs{j}}\Gless\cGs{\xs{\is{k+1}}}}$.
		
		We claim the subsequence of points chosen this way still has the property that successive points 
		belong to a common strip: if $\is{k+1}>\is{k}+1$, then $\xs{\is{k+1}-1}\leq\xs{\is{k}}\leq\xs{\is{k+1}}$,
		so the strip containing \xs{\is{k+1}-1} and \xs{\is{k+1}} contains \xs{\is{k}}.
		
		But now the \Gless-monotonicity of the sequence means that the transversals \clint{\xs{\is{k-1}}}{\xs{\is{k}}}
		and \clint{\xs{\is{k}}}{\xs{\is{k+1}}} lie on opposite sides of \cGs{\xs{\is{k}}}, and hence their union is still
		transverse to \cG{}.  But then it is a global cross section for \U, establishing \refer{lem}{connect} and hence
		\refer{prop}{sep}.		   

\end{proofof}

\begin{proofof}{Theorem C}
	To prove Theorem C, we need to show that if \Realstwo{} does not have the separation property for \cG{}
	then it contains a Reeb component in the sense of \refer{dfn}{Reeb}.
	
	Consider the collection of open, connected \cG-saturated sets with the separation property:  it is clearly nonempty,
	and any nested union of such sets is again such a set.  Therefore by Zorn's lemma we can find a \emph{maximal}
	open, connected \cG-saturated set \U{} with the separation property.  By assumption, $\U\neq\Realstwo$,
	and hence has a nonempty boundary $\bdU$. Note that this boundary is \cG-saturated: for each $p\in\bdU$,
	$\cGs{p}\subset\bdU$.
	
	We claim that there exists a pair of distinct leaves $\cGs{x},\cGs{y}\subset\U$ such that $p$ and $x$ 
	\resp{$p$ and $y$} lie on the same side of \cGs{y} \resp{\cGs{x}}--because if not, then since a tubular
	neighborhood of \cGs{p} has the separation property and $p$ separates any point outside \U{} from any point
	inside \U, the union of \U{} with the interior of a tubular neighborhood of \cGs{p} would contradict the maximality of \U.
	
	So fix two such leaves $\cGs{x}\neq\cGs{y}\subset\U$;  we can assume $\cGs{x}\Gless\cGs{y}$.
	
	For each point $z\in\U$ with $\cGs{x}\Gless\cGs{z}\Gless\cGs{y}$, the complementary planes satisfy
	\begin{align*}
		\sideneg{x}\subset&\sideneg{z}\subset\sideneg{y}\\
		\sidepos{x}\supset&\sidepos{z}\supset\sidepos{y}.
	\end{align*}
	Note that 
	\begin{itemize}
		\item If $p\in\sidepos{z}$, then \cGs{z} separates $p$ from \cGs{x}, while
		\item If $p\in\sideneg{z}$, then \cGs{z} separates $p$ from \cGs{y}.
	\end{itemize}
	Let $X$ be the union for all $z$ in the first situation of the \emph{negative} sides \sideneg{z} 
	and $Y$ the union of  \emph{positive} sides \sidepos{z} for $z$ in the second situation (note that $x\in X$ and $y\in Y$).
	
	Consider a transversal arc $\gamma=\clint{x}{y}$; 
	this is a global cross section for the region $\U\cap\sidepos{x}\cap\sideneg{y}$, 
	whose boundary includes \cGs{p}, \cGs{x} and \cGs{y};
	note that \gam{} is oriented positively (in terms of \Gless) from $x$ to $y$.
	
	The two intersections $X\cap\gamma$ and $Y\cap\gamma$ are 
	connected and complementary, so there is a dividing point 
	\begin{equation*}
		z=\sup(X\cap\gam)=\inf(Y\cap\gam)\in\gamma
	\end{equation*} 
	which belongs to one of them, but is a 
	limit point of the other.  Suppose for definiteness that $z\in Y$, so \cGs{z} separates $p$ from \cGs{y}.
		
	Now take a sequence \xs{k} of points in \U{} converging to $p$; they must belong to $X$, and so $p\in\sidepos{\xs{k}}$.
	This means that we can assume without loss of generality that the sequence \cGs{\xs{k}} is \cG-increasing.  
	Now consider the intersections $\xps{k}=\cGs{\xs{k}}\cap\gamma$. These are also
	\Gless-increasing, and have a supremum.  It cannot fall short of $z$, since then the leaves crossing \gam{} between
	this supremum and $z$ separate $p$ from \xs{k}.
	
	So we have a sequence of pairs  of points,$(\xs{k}, \xps{k})$, with $\xs{k}\to p$ and $\xps{k}\to z$, 
	but $p$ and $z$ lie on different leaves;  hence they are prolongationally related, so the region
	bounded by \cGs{p} and \cGs{z} is a Reeb component in the sense of \refer{dfn}{Reeb}.
	This proves Theorem C.
	
\end{proofof}
\subsection{Action of a time-one map on a foliation}\label{subsec:action}
We now consider the following situation:
\begin{itemize}
	\item \flowt{} is a fixedpoint-free flow on \Realstwo;
	\item \cG{} is the foliation of \Realstwo{} by orbits of the flow;
	\item $f=\flowto{1}$ is the time-one map of the flow;
	\item \cF{} is a foliation preserved by $f$.
\end{itemize}
We emphasize that \cF{} is \emph{not} assumed to be preserved by the flow \flowt{} for non-integer times $t$%
\footnote{
In fact, it was the distinction between $f$-preservation and \flowt-preservation that allowed White to construct his example of an
Anosov structure for a translation.
}.

The main theorem of this section is
\begin{ThmB}
	Suppose \flowt{} is a fixedpoint-free \Cr{1} flow in the plane with flow line foliation \cG.  
	Let \selfmap{f=\flowto{1}}{\Realstwo} be the time-one map
	of \flowt, and suppose \cF{} is a \Cr{1} foliation preserved by $f$.
	
	If \cG{} is non-trivial, then some leaf of \cF{} is invariant under $f$.
\end{ThmB}

%
\subsubsection{Relation between leaves of two different foliations}
Before proving Theorem B, we consider some technical results involving an interplay of leaves for two different 
\Cr{1} foliations.  
We  apply all the notation developed in \refer{subsec}{PlanFolPrelims} to both foliations \cG{} and \cF{}, distinguishing which
foliation is involved via the subscript.

Suppose that $p\in\Realstwo$ is a point with $\cFs{p}\neq\cGs{p}$;  by picking a point $\xp\in\cFs{p}\setminus\cGs{p}$ 
and then taking a maximal subarc \Fropint{x}{\xp}  disjoint from \cGs{p}, we have a point $x\in\cGs{p}$ and two arcs,
\Fclint{x}{\xp} and \Gclint{x}{y} such that 
\begin{equation*}
	\Fclint{x}{\xp}\cap\cGs{p}=\single{x}=\Gclint{x}{y}\cap\cFs{p}.
\end{equation*}
%
%


\begin{prop}\label{prop:crossing}
	Given two arcs \Fclint{x}{\xp} and \Gclint{x}{y} satisfying the condition above, any neighborhood $U$ of $x$ contains a 
	(topological) disc $D$ separated by a \Garc{} $I\subset\Gclint{x}{y}$, such that every point $q\in D$ is joined to
	a point of $I$ by an \Farc{} in $D$.
\end{prop}

\begin{proofof}{\refer{prop}{crossing}}

We can assume that \Fclint{x}{\xp} is contained in $U$.  Given a sufficiently narrow \cF-half-tubular neighborhood of \cFs{p},
not including $y$,
we can take an arc \clint{\xp}{\xps{1}} in this tubular neighborhood transverse to \cF{}. 
Also, letting \vvs{\cF} \resp{\vvs{\cG}} be the unit vector tangent at $x$ to \Fclint{x}{\xp} \resp{\Gclint{x}{y}},
we can pick a vector \vvs{1} in the  sector bounded by \vvs{\cG} and $-\vvs{\cF}$;  
since both foliations are \Cr{1}, a sufficiently short straight line segment \clint{x}{\xs{1}}from $x$ with direction vector \vvs{1} constitutes
an arc \clint{x}{\xs{1}} transverse to 
both \cF{} and \cG{}, contained in the component of the complement of \cGs{p} not containing \xp.
Let $D$ be the \Fbox{} cut out from the \cF-half-tubular neighborhood by \clint{\xp}{\xps{1}} and \clint{x}{\xs{1}},
with \Fclint{\xs{1}}{\xps{1}} redefined to be the edge of $D$ opposite \Fclint{x}{\xp}.  
Then set $I$ to be the component of $\cGs{p}\cap D$ containing $x$.  This is a subarc of \clint{x}{y}, so does not intersect \Fclint{x}{\xp}
or either of the transversals;  but $y$ does not belong to $D$, so the other endpoint of $I$ (besides $x$) must lie on \Fclint{\xs{1}}{\xps{1}},
and hence $I$ separates $D$.  Since the two transversals lie on opposite sides of $I$ and every \cF-leaf intersects $D$ in a path joining the
two transversals, every point in $D$ is joined to $I$ by an \Farc{} inside $D$, as required.
\end{proofof}

\subsubsection{Action}
\refer{dfn}{proldyn} can be repeated for the discrete dynamical system generated by a homeomorphism $f$ of \Realstwo{}
by replacing the ``times'' $\tsi\in\Reals$ with ``iterates'' $\ksi\in\Integers$.  When $f=\flowto{1}$ is the time-one map of \cG, 
we need to be 
careful to distinguish prolongation under \flowt{} from prolongation under $f$.  In general, $\prolpf{x}\subset\prolpfl{x}$,
but the reverse inclusion does not hold in general..  However, we have the following

\begin{prop}\label{prop:prolftvsf}
	Suppose $f=\flowto{1}$ where \flowt{} is a fixedpoint-free flow in \Realstwo, with flowline foliation \cG.
	If $y\in\prolpfl{x}$, then \prolpf{x} intersects the \Garc{} \flowclint{y}{y+1}.
\end{prop}
%

\refer{prop}{prolftvsf} will follow from the following apparently weaker result.
\begin{lemma}\label{lem:flvf1}
	Given $y\in\prolpfl{x}$, pick \epsgo{} and consider the closed \cG-arc 
	\begin{equation*}
		\Is{\eps}\eqdef\flowclint{y-\eps}{y+1+\eps}.
	\end{equation*}
	Then \prolpf{x} intersects \Is{\eps}.
\end{lemma}

\begin{subproof}{\refer{lem}{flvf1}}

	Note that any \cG-arc in \cGs{y} of \Glen{} at least $1$ contains a point of the $f$-orbit of any point in \cGs{y}.

	Let \Ts{0} \resp{\Ts{1}} be arcs through $\flowto{-\eps}(y)$ \resp{$\flowto{1+\eps}(y)$}, transverse to \cG.
	By \refer{rmk}{Gbox}, shrinking these transversals if necessary, 
	we can assume they are the vertical sides of a \Gbox{} containing
	\Is{\eps}.  Since the \Glen{} of horizontal arcs in a \Gbox{} varies continuously and equals $1+2\eps$
	for \Is{\eps}, if the height of this \Gbox{} is sufficiently small,
	its inner \Glen{} is at least $1+\eps$..
	
	In particular, since $y\in\prolpfl{x}$, there are points \xs{i} converging to $x$ (in \U) and points 
	$\ysi\in\cGs{\xsi}$ converging to $y$.  This means the horizontal \Garc{}s through \ysi{} converge to
	\Is{\eps}.  But each of these arcs contains a point \zsi{} in the $f$-orbit of \xsi{}, and by compactness these
	have at least one accumulation point $z$ which then belongs to $\Is{\eps}\cap\prolpf{x}$.
\end{subproof}

\begin{proofof}{\refer{prop}{prolftvsf}}
	Since \prolpf{x} is closed, so is its intersection with each \Garc{} \Is{\eps} for $\eps\to0$.
	Hence  by the nested intersection property
	the intersection with $\bigcap_{\eps}\Is{\eps}=\flowclint{y}{y+1}$ is nonempty.
\end{proofof}

%
%
%
%
%
%
\begin{lemma}\label{lem:freeleaf}
	An arc joining the edges of a Reeb component of \cG{} must intersect its $f$-image.
\end{lemma}

\begin{subproof}{\refer{lem}{freeleaf}}
	Suppose for some Reeb component of \cG, $x\in\Gamn$ and $\xp\in \Gamp$, and \gam{} is an arc with endpoints
	$\xs{-}\in\Gamn$ and $\xs{+}\in\Gamp$.  Since \Gamn{} and \Gamp{} are closed and disjoint, 
	we can assume that the interior of \gam{} is contained
	in \U, the component of the complement of $\Gamn\cap\Gamp$ bounded by these two curves.  
	Let \gams{\pm} be the component of $\Gams{\pm}\setminus\{\xs{\pm}\}$ containing \fof{\xs{\pm}}.
	Then $\gams{0}\eqdef\gams{-}\cup\gam\cup\gams{+}$ separates \U, and if $\fof{\gam}\cap\gam=\emptyset$
	then the component of $\U\setminus\gamma$ containing \fof{\xs{-}} is mapped into itself under $f$.  Then
	\prolpf{\xs{-}} cannot include any point in the other component of $\U\setminus\gamma$, contradicting 
	\refer{prop}{prolftvsf} with $x=\xs{-}$ and $y=\ftoof{-1}{\xs{+}}$.
\end{subproof}

\begin{proofof}{Theorem B} 
Clearly, we can assume that no leaf of \cF{} is also a leaf of \cG{}, since each leaf of \cG{} is $f$-invariant. 
So suppose $f$ is the time-one map of a fixedpoint-free flow \Cr{1} \flowt{} whose flow line foliation \cG{} is nontrivial, 
and hence by Theorem C has a Reeb component, formed by \Gamn, \U, and \Gamp{} as in \refer{dfn}{Reeb}.
Suppose furthermore that \cF{} is another
\Cr{1} foliation with no $f$-invariant leaves--in particular, no leaf of \cF{} coincides with any leaf of \cG.

Since the \cF-leaf of any point 
of \Gamn{} is not equal to \Gamn{}, by \refer{prop}{crossing} there exist an open disc \Dm{} in \Realstwo{} intersecting \Gamn{} 
and a compact interval $I\subset\Gamn$ such that
any \cF-leaf intersecting $\Dm{}\setminus I$ contains an arc, intersecting $I$, whose interior is disjoint from \Gamn.  

Pick $x\in\Dm\cap\Gamn$ and $y\in\prolpf{x}$.  
Note that $y\in\prolpfl{x}\subset\Gamp$.  
Let \Ts{0} be an arc in the closure of \U{}, transverse to both \cG{} and \cF{},
with one endpoint at $y$ and the rest internal to \U.  Since $f$ is a diffeomorphism, $\Ts{1}=\fof{\Ts{0}}$ is also transverse 
to both foliations, has an endpoint at \fof{y}, and the rest of it is contained in the $f$-invariant set \U.

Since $y\in\prolpf{x}$, there exist points \xsi{} in $(\Dm{}\setminus I)\cap\U$, converging to $x$, and times $\ksi\to+\infty$ 
such that $\ysi\eqdef\ftoof{\ksi}{\xsi}\to y$.  Since \U{} is $f$-invariant and separates \Realstwo, \xsi{} and \ysi{} belong to \U.
Since $\xsi\in\Dm\setminus I$, 
there is a point $\xpsi\in I$ 
such that \Fopint{\xsi}{\xpsi} is contained in \U.

Note that, since $f$ restricts to a fixedpoint-free homeomorphism of \Gamn,
there is a positive iterate $f^{k}$ such that the $f^{k}$-images of $I$ are pairwise disjoint.  
Furthermore, since $f$ preserves order in \Gamn,
if $\ksi<\ks{i+1}<\ks{i+2}$ with $\ks{i+1}-\ksi>k$ and $\ks{i+2}-\ks{i+1}>k$ then setwise \ftoof{\ks{i+1}}{I} is between \ftoof{\ksi}{I}
and \ftoof{\ks{i+2}}{I}. 

Since \ysi{} are eventually in any \Fbox{}
neighborhood of $y$ and \Ts{0} is transverse to \cF, we can assume  (passing to a subsequence) that each \cF-leaf
\cFs{\ysi} intersects \Ts{0} at a point \ypsi. 
Letting $\zsi=\ftoof{\ksi}{\xpsi}$, we see that 
$\Fopint{\ypsi}{\zsi}=\ftoof{\ksi}{\Fopint{\xsi}{\xpsi}}$ is contained in \U.

\begin{lemma}\label{lem:ypsimono}
	The sequence $\{\ypsi\}$ converges to $y$ monotonically in \Ts{0}.
\end{lemma}
\begin{subproof}{\refer{lem}{ypsimono}(see \refer{fig}{ypsimono})}
Passing to a further subsequence, we assume that for every $i$, $\ks{i+1}-\ksi>k+2$,  
so that the points \zsi{} diverge monotonically to $+\infty$ in the orientation of \Gamn{}
induced by \flowt, and for each $i$, \fof{\zsi} and \ftoof{2}{\zsi} lie between \zsi{} and \zs{i+1}.

Now fix $i$ momentarily, and consider the curve \gam{} consisting of 
$\gams{i}=\Fclint{\ypsi}{\zsi}$, 
the arc in \Ts{0} with endpoints at \ypsi{} and $y$, 
and the arcs \Gropint{\zsi}{+\infty} and \Glopint{-\infty}{y} in \Gamn{} and \Gamp{}, respectively%
\footnote{
Here we are using the natural ordering of \Gamn{} and \Gamp{} induced by the flow to define ``$\pm\infty$'' in either leaf.
}.
\gam{} separates \Realstwo; let $V$ be the component of $\U\setminus\gam$ containing \fof{\zsi} (and hence also \ftoof{j}{\zsi}
for all $j\geq1$).
Note that \Ts{1} lies outside $V$, so $\fof{\gams{i}}=\Fopint{\fof{\ypsi}}{\fof{\zsi}}$ must cross \gam. It can't intersect either
\gams{i} (because this is part of a different leaf of \cF) or \Gams{\pm} (since it is contained in the $f$-invariant set \U).
Thus, it must cross \Ts{0} between \ypsi{} and $y$, say at $p(i,1)$.  Then any \Farc{} joining a point of \Glopint{\fof{\zsi}}{+\infty}
to \Ts{1} must intersect \Ts{0} between $p(i,1)$ and $y$.
In particular this is true of $\ftoof{j}{\gams{i}}=\Fopint{\ftoof{j}{p(i,1)}}{\ftoof{j}{\zsi}}$ for every $j>0$
as well as $\gams{i+1}=\Fopint{\yps{i+1}}{\zs{i+1}}$ and its $f^{j}$-images.


An analogous argument shows that if three \Farc{s} \Fclint{\asi}{\bsi} $i=1,2,3$ in \U{} have $\asi\in\Gamn$ and $\bsi\in\Ts{0}$
and \as{2} is between \as{1} and \as{3} in \Gamn, then \bs{2} is between \bs{1} and \bs{3} in \Ts{0}. 

For $i=1,2,...$ and $j=0,1,2,...$, let $p(i,j)$ be the intersection of \ftoof{j}{\gams{i}} with \Ts{0} (of course $p(i,0)=\ypsi$). 
For each fixed $i$, the sequence \ftoof{j}{\zs{i}} is monotone increasing in \Gamn, so the sequence 
$p(i,1),p(i,2),...$ is monotone in \Ts{0}. 
Similarly, since the sequence \zs{i} is monotone increasing in \Gamn, for fixed $j$ the sequence $p(1,j),p(2,j),...$ is monotone in \Ts{0}.  Finally, since $\ys{i}\to y$ eventually belong to an \Fbox{} around $y$, $p(i,0)\to y$ as $i\to\infty$, so
$\fof{p(i,0)}\to\fof{y}$. \end{subproof}

%
%

\begin{figure}[htbp]
\begin{center}
	\begin{pspicture}(-4,-2.5)(8.5,4)
		\pscircleOA*[linecolor=gray!20](-3,3)(-4,3)
			\uput[d](-2.4,3){$\boldsymbol{x}$}\rput(-2.4,3){$\bullet$}
			\rput(-3.5,2.2){\xs{i}} \uput[u](-3.5,3){\xps{i}} \psline[linecolor= gray](-3.5,2.2)(-3.5,3)
			\rput(-2.8,2.5){\xs{i+1}}\uput[u](-2.8,3){\xps{i+1}}\psline[linecolor= gray](-2.8,2.5)(-2.8,3)

		\psline[linewidth=1.5pt]{->}(-4,3)(8,3)
			\uput[r](8,3){$\Gamma_{-}	$}
		\psline[linewidth=1.5pt]{<-}(-4,-2)(8,-2)
			\uput[r](8,-2){$\Gamma_{+}$}

		\uput[dr](0,-2){$y\in J{f}^{+}(x)$}\rput(0,-2){$\bullet$}
			\psline(0,-2)(0,1)
			\uput[l](0,0.5){\ys{i}}
			\uput[r](0,0.5){\ftoof{\ks{i}}{\xs{i}}}
		
		\uput[d](-1.5,-2){$f(y)$}
			\psline(-1.5,-2)(-1.5,1)
			\uput[l](-1.5,-0.3){\fof{\ys{i}}}
			\uput[r](-1.5,-0.3){\ftoof{\ks{i}+1}{\xs{i}}}
		
		\uput[d](-3,-2){$f^{2}(y)$}
			\psline(-3,-2)(-3,0)
			\uput[l](-3,-1.0){\ftoof{2}{\ys{i}}}
			\uput[r](-3,-1.0){\ftoof{\ks{i}+1}{\xs{i}}}
		
		\uput[u](1.5,3){\ftoof{\ks{i}}{\xps{i}}} 
			\pscurve[linecolor=gray](1.5,3)(0.5,0.2)(0,0.2)

		\uput[u](3.5,3){\ftoof{\ks{i}+1}{\xps{i}}} 
			\pscurve[linecolor=gray](3.5,3)(0,-0.6)(-1.0,-0.6)(-1.5,-0.6)

		\uput[u](5.5,3){\ftoof{\ks{i}+2}{\xps{i}}} 
			\pscurve[linecolor=gray](5.5,3)(0,-1.3)(-0.7,-1.0)(-1.5,-1.3)(-2.5,-1.3)(-3,-1.3)

		\uput[u](7.5,3){\ftoof{\ks{i+1}}{\xps{i+1}}} 
			\pscurve[linecolor=gray](7.5,3)(0,-1.5)(-0.7,-1.3)(-1.5,-1.5)(-2.5,-1.5)(-3,-1.7)(-3.5,-1.7)

			\pscurve[linecolor=black, linestyle=dashed, linewidth=1.2pt](8,2.5)(0,-1.9)(-0.7,-1.85)(-1.5,-1.9)(-2,-1.7)(-3,-1.9)(-3.5,-1.8)

	\end{pspicture}
\caption{Proof of Lemma~\ref{lem:ypsimono}}
\label{fig:ypsimono}
\end{center}
\end{figure}

To complete the proof of Theorem B, note that all the \Farc{}s $\alps{i}=\Fclint{\fof{p(i,0}}{p(i,1)}$ 
are contained in the topological closed rectangle \FR{} bounded
by \alps{1}, \Fclint{y}{\fof{y}}, and subintervals of \Ts{0} and \Ts{1}.  
All the \alps{i} are disjoint, and their endpoints in each of \Ts{0} and \Ts{1}
converge monotonically to $y$ and \fof{y}, respectively; it follows that \alps{i} converge in the Hausdorff topology 
to an \Farc{} \alps{\infty} whose endpoints are $y$ and \fof{y}.  But this says that $\cFs{y}=\cFs{\fof{y}}$, making it an 
$f$-invariant leaf of \cF.
\end{proofof}

%

\section{The Mendes Conjecture-a partial resolution}\label{sec:Mendes}

In this section we prove Theorem A.
We separate the two cases: (1) $f$ is fixedpoint free, and (2) $f$ has a unique nonwandering point. 

\subsection{Case 1: $f$ is a Brouwer translation}
When \selfmap{f}{\Realstwo} is fixedpoint free (a Brouwer translation) then the Mendes conjecture says that $f$ must be
topologically conjugate to the translation%
\footnote{
Note that all translations are mutually topologically conjugate.
} 
$(x,y)\mapsto(x+1,y)$.
Under the additional assumption that $f$ is the time-one map of some flow, using Theorem B we have
the following

\begin{prop}\label{prop:trivial}
	If \flowt{} is a flow with trivial orbit foliation then there is a homeomorphism \selfmap{h}{\Realstwo} which is equivariant
	with respect to \flowt{} and the translation flow $\Phi^{t}$ defined by $\Phi^{t}(x,y)=(x+t,y)$:
	\begin{equation}\label{eqn:flowconj}
		\hof{\flowtoof{t}{(x,y)}}=\Phi^{t}(\hof{x,y})\text{ for all $(x,y)\in\Realstwo$ and $t\in\Reals$}.
	\end{equation}
\end{prop}

\begin{proofof}{\refer{prop}{trivial}}
	For flows in the plane, triviality of the orbit foliation of \flowt{} is equivalent to the existence of a single, connected 
	cross-section--a line meeting every orbit at exactly once.  Pick such a section $T$ for \flowt{} 
	and initially define $h$ on $T$ to be any homeomorphism between $T$ and the \axis{y} in \Realstwo. 
	Then we extend the definition of $h$ to the whole plane by noting that for each $(x,y)\in\Realstwo$ 
	there is a (unique) point $(\xp,\yp)\in T$ and $t\in\Reals$ such that $\varphi^{t}(\xp,\yp)=(x,y)$; by definition, we want
	\begin{equation*}
		\hof{(x,y)}=\Phi^{t}(\hof{\xp,\yp})
	\end{equation*}
	which gives the required conjugacy.
\end{proofof}

\begin{corollary}\label{cor:Mendes1}
	If the time-one map $f$ of a fixedpoint-free flow \flowt{} in \Realstwo{} has an Anosov structure, then the 
	action of the flow is conjugate to that of the translation flow, and so $f$ is topologically conjugate
	to the translation $T$.
\end{corollary}

\begin{proofof}{\refer{cor}{Mendes1}}
	Suppose the flowline foliation of \flowt{} is non-trivial, and $f$ has an Anosov structure.
	Let \cF{} be the associated stable foliation of \Realstwo{}.
	Clearly, \cF{} is $f$-invarient, so by Theorem B some leaf of \cF{} is $f$-invariant.
	But then $f$ restricted to this leaf is a contraction with respect to the metric giving the Anosov structure,
	and hence has a fixedpoint, contrary to the assumption that \flowt{} is fixedpoint-free.
	Thus, \cG{} must be trivial. But then by \refer{prop}{trivial} there is a homeomorphism \selfmap{h}{\Realstwo}
	 such that \refer{eqn}{flowconj} holds;  in particular, setting $t=1$, we get the conjugacy condition
	\begin{equation*}
		\hof{\fof{(x,y)}}=\Phi^1(\hof{x,y})\text{ for all $(x,y)\in\Realstwo$}.
	\end{equation*}
\end{proofof}

\subsection{Case 2: $f$ has a fixedpoint}\label{subsec:fixedpoint}

Our standing assumption in this subsection is that $f$ is the time-one map of a \Cr{1} flow \flowt{} on \Realstwo,
has an Anosov structure, 
and has a unique fixedpoint.  


The second condition in \refer{thm}{MendesThm} implies 
that the stable and unstable manifolds of this fixedpoint escape to infinity.  
Then  the ``cross'' $X$ consisting of the fixedpoint and its separatrices 
separates the plane into four $f$-invariant open quadrants $\Qsi$, $i=1,...4$.  
If $f$ is the time-one map of a flow \flowt, 
the restriction of the flow to each quadrant is fixedpoint-free,  hence generates a flowline foliation \cG.  
By Theorem B, if \cG{} is nontrivial, then any foliation in this quadrant which is preserved by $f$ 
must have an $f$-invariant leaf.  Applying this to the foliation (of the open quadrant) by stable manifolds, we would have to have
a second fixedpoint of $f$, contrary to the first condition in \refer{thm}{MendesThm}. It follows that the foliation \cG{} restricted to each quadrant must be trivial: 
\begin{remark}\label{rmk:quad}
	There is a homeomorphism of each open quadrant \Qsi{} to \Realstwo{} 
	(which here we represent as the open upper half-plane) 
	taking flow lines of \flowt{} to horizontal lines.
\end{remark}

We will construct a conjugating homeomorphism \selfmap{h}{\Realstwo} using a standard trick.  
We call a closed topological disc $D\subset\Realstwo$ a \deffont{fundamental domain} for a homeomorphism
\selfmap{g}{\Realstwo} if there are two closed arcs \gams{-} and \gams{+} in its boundary such that
\begin{equation*}
	\gams{+}=\gof{\gams{-}}=D\cap\gof{D}.
\end{equation*} 
\begin{remark}\label{rmk:fundom}
	If $D$ \resp{$\tilde{D}$} is a fundamental domain for $g$ \resp{\tg} and \map{h}{D}{\tilde{D}} is a homeomorphism
	taking \gams{\pm} to $\tilde{\gamma}_{\pm}$, then $h$ extends to a homeomorphism conjugating
	$g|\bigcup_{k\in\Integers}\toof{g}{k}{D}$ with $\tg|\bigcup_{k\in\Integers}\toof{\tg}{k}{\tilde{D}}$ via
	\begin{equation*}
		\hof{\toof{g}{k}{x}}=\toof{\tg}{k}{\hof{x}} \text{ for all } x\in D.
	\end{equation*}
\end{remark}

We let \cFst{} \resp{\cFut} be the foliation of \Realstwo{} by the stable \resp{unstable} manifolds of $f$. Note that the stable
\resp{unstable} separatrices of the fixedpoint are contained in a leaf of  \cFst{} \resp{\cFut}.
\begin{lemma}\label{lem:fbox}
	There is a rectangular neighborhood \cR{} of the fixedpoint of $f$
	whose horizontal \resp{vertical} edges are \cFut-arcs \resp{\cFst-arcs}, 
	which is simultaneously an \cFst-box and a \cFut-box. 
\end{lemma}
\begin{subproof}{\refer{lem}{fbox}}
	Since the \cFst-leaf and \cFut-leaf through the fixedpoint are transversal (and consist, respectively, of the appropriate
	separatrices together with the fixedpoint itself), there is a disc neighborhood of the fixedpoint  on which the two
	foliations form a product structure: any \cFst-arc in the neighborhood intersects any \cFut-arc in at most one point
	(and is transversal).  Pick a pair of \cFut-arcs, one through a point on each stable separatrix.
	\cFst-arcs of points near the fixedpoint intersect both arcs; pick one such  \cFst-arc through a point on each unstable
	separatrix.  There are four points \qs{i} of intersection between the two \cFst-arcs and the two \cFut-arcs; with 
	appropriate numbering the rectangle \cR{} formed by 
	$\Eua=\Fuclint{\qs{1}}{\qs{2}}$, $\Esa=\Fsclint{\qs{2}}{\qs{3}}$, 
	$\Eub=\Fuclint{\qs{3}}{\qs{4}}$,
	and $\Esb=\Fsclint{\qs{4}}{\qs{1}}$ is foliated by \cFut-arcs joining the two ``vertical'' edges \Esa{} and
	\Esb, and also by the \cFst-arcs joining the ``horizontal'' edges \Eua{} and \Eub. 	
\end{subproof}

For future reference, we note that each of these four edges crosses one of the separatrices of the fixedpoint at a unique point;
denote the ``cross'' formed by the fixedpoint together with its four separatrices by $X$ and set
\begin{align*}
	\ps{1}&=X\cap\Eua\\
	\ps{2}&=X\cap\Esa\\
	\ps{3}&=X\cap\Eub\\
	\ps{4}&=X\cap\Esb.
\end{align*}

\begin{figure}[htbp]
\begin{center}
	\begin{pspicture}(-4,-4)(4,4)
\psset{unit=2}
\pscustom[linestyle=dashed, linecolor=black, fillstyle=solid, fillcolor=gray!20]{\pspolygon%
		(-1.0,-1.0)(-0.5,-0.9)(-0.2,-1.0)(0,-1.0)(0.2,-1.0)(0.5,-0.9)(1.0,-1.0)%
		(1.0,-1)(1.1,-0.5)(1.0,-0.2)(1.0,0)(1.0,0.2)(1.1,0.5)(1.0,1)%
		(0.5,1.1)(0.2,1.0)(0,1.0)(-0.2,1.0)(-0.5,1.1)(-1.0,1.0)%
		(-0.9,0.5)(-1.0,0.2)(-1.0,0)(-1.0,-0.2)(-0.9,-0.5)(-1.0,-1)
		}

		\rput(0,0){$\bullet$}

		\pscurve[linecolor=black, linewidth=1.5pt,arrows=<->, arrowscale=1.5](-2,0)(-1.5,-0.1)(-1,0)(-0.5,0.1)(-0.2,0)(0,0)(0.2,0)(0.5,0.1)(1,0)(1.5,-0.1)(2,0)

		\pscurve[linecolor=gray, linewidth=1.5pt,arrows=>-<, arrowscale= 1.5](0,-2) (-0.1,-1.5)(0,-1)(0.1,-0.5)(0,-0.2)(0,0)(0,0.2) (0.1,0.5)(0,1)(-0.1,1.5)(0,2)

	\pscurve[linecolor=gray](0.2,-1)(0.3,-0.5)(0.2,-0.2)(0.2,0)(0.2,0.2)(0.3,0.5)(0.2,1)\uput[d](0.2,-1){\ps{3}}\uput[u](0.2,1){\ps{1}}
	\pscurve[linecolor=gray](0.4,-1)(0.5,-0.5)(0.4,-0.2)(0.4,0)(0.4,0.2)(0.5,0.5)(0.4,1)
	\pscurve[linecolor=gray](0.6,-1)(0.7,-0.5)(0.6,-0.2)(0.6,0)(0.6,0.2)(0.7,0.5)(0.6,1)
	\pscurve[linecolor=gray](0.8,-1)(0.9,-0.5)(0.8,-0.2)(0.8,0)(0.8,0.2)(0.9,0.5)(0.8,1)
	\pscurve[linecolor=gray](1.0,-1)(1.1,-0.5)(1.0,-0.2)(1.0,0)(1.0,0.2)(1.1,0.5)(1.0,1)\uput[r](1.0,-1){\qs{3}}\uput[r](1.0,1){\qs{2}}

	\pscurve[linecolor=gray](-0.2,-1)(-0.1,-0.5)(-0.2,-0.2)(-0.2,0)(-0.2,0.2)(-0.1,0.5)(-0.2,1)
	\pscurve[linecolor=gray](-0.4,-1)(-0.3,-0.5)(-0.4,-0.2)(-0.4,0)(-0.4,0.2)(-0.3,0.5)(-0.4,1)
	\pscurve[linecolor=gray](-0.6,-1)(-0.5,-0.5)(-0.6,-0.2)(-0.6,0)(-0.6,0.2)(-0.5,0.5)(-0.6,1)
	\pscurve[linecolor=gray](-0.8,-1)(-0.7,-0.5)(-0.8,-0.2)(-0.8,0)(-0.8,0.2)(-0.7,0.5)(-0.8,1)
	\pscurve[linecolor=gray](-1.0,-1)(-0.9,-0.5)(-1.0,-0.2)(-1.0,0)(-1.0,0.2)(-0.9,0.5)(-1.0,1)\uput[l](-1,-1){\qs{4}}\uput[l](-1,1){\qs{1}}

	\pscurve[linecolor=black](-1.0,0.2)(-0.5,0.3)(-0.2,0.2)(0,0.2)(0.2,0.2)(0.5,0.3)(1.0,0.2)
	\pscurve[linecolor=black](-1.0,0.4)(-0.5,0.5)(-0.2,0.4)(0,0.4)(0.2,0.4)(0.5,0.5)(1.0,0.4)
	\pscurve[linecolor=black](-1.0,0.6)(-0.5,0.7)(-0.2,0.6)(0,0.6)(0.2,0.6)(0.5,0.7)(1.0,0.6)
	\pscurve[linecolor=black](-1.0,0.8)(-0.5,0.9)(-0.2,0.8)(0,0.8)(0.2,0.8)(0.5,0.9)(1.0,0.8)
	\pscurve[linecolor=black](-1.0,1.0)(-0.5,1.1)(-0.2,1.0)(0,1.0)(0.2,1.0)(0.5,1.1)(1.0,1.0)
	
	\uput[ul](-1,0){\ps{4}}\uput[ur](1,0){\ps{2}}

	\pscurve[linecolor=black](-1.0,-0.2)(-0.5,-0.1)(-0.2,-0.2)(0,-0.2)(0.2,-0.2)(0.5,-0.1)(1.0,-0.2)
	\pscurve[linecolor=black](-1.0,-0.4)(-0.5,-0.3)(-0.2,-0.4)(0,-0.4)(0.2,-0.4)(0.5,-0.3)(1.0,-0.4)
	\pscurve[linecolor=black](-1.0,-0.6)(-0.5,-0.5)(-0.2,-0.6)(0,-0.6)(0.2,-0.6)(0.5,-0.5)(1.0,-0.6)
	\pscurve[linecolor=black](-1.0,-0.8)(-0.5,-0.7)(-0.2,-0.8)(0,-0.8)(0.2,-0.8)(0.5,-0.7)(1.0,-0.8)
	\pscurve[linecolor=black](-1.0,-1.0)(-0.5,-0.9)(-0.2,-1.0)(0,-1.0)(0.2,-1.0)(0.5,-0.9)(1.0,-1.0)

	\rput(0.5,1.5){$\boldsymbol{\Eua}$}\psline[linestyle=dotted](0.5,1.5)(0.5,1.1)
	\rput(1.5,0.5){$\boldsymbol{\Esa}$}\psline[linestyle=dotted](1.5,0.5)(1.0,0.4)
	\rput(-0.5,-1.5){$\boldsymbol{\Eub}$}\psline[linestyle=dotted](-0.5,-1.5)(-0.5,-1.1)
	\rput(-1.5,-0.5){$\boldsymbol{\Esb}$}\psline[linestyle=dotted](-1.5,-0.5)(-1.0,-0.4)

	\end{pspicture}
\caption{The rectangle \cR}
\label{fig:rect}
\end{center}
\end{figure}

We now form a larger neighborhood \Vso{} of the fixedpoint by first taking the union $\cR\cup\fof{\cR}$, then further enlarging
by joining each vertex \qsi{} of \cR{} with its image \fof{\qsi} by the \cG-arc $\cGs{i}=\Gclint{\qsi}{\fof{\qsi}}$;
the resulting topological octagon 
\begin{equation*}
	\partial \Vso=\Eua\cup\cGs{1}\cup\fof{\Esa}\cup\cGs{2}\cup\Eub\cup\cGs{3}\cup\fof{\Esb}\cup\cGs{4}
\end{equation*} 
bounds a closed topological disc \Vso{} which is also simultaneously an \cFst-box and a \cFut-box 
(provided our initial choices were sufficiently close to the fixedpoint).

\begin{figure}[htbp]
\begin{center}
	\begin{pspicture}(-4,-4)(4,4)
\psset{unit=2}
\pscustom[linestyle=none, fillstyle=solid, fillcolor=gray!20]%
	{\pspolygon%
		(-1.0,-1.0)(-0.5,-0.9)(-0.2,-1.0)(0,-1.0)(0.2,-1.0)(0.5,-0.9)(1.0,-1.0)%
		(1.0,-1)(1.1,-0.5)(1.0,-0.2)(1.0,0)(1.0,0.2)(1.1,0.5)(1.0,1)%
		(0.5,1.1)(0.2,1.0)(0,1.0)(-0.2,1.0)(-0.5,1.1)(-1.0,1.0)%
		(-0.9,0.5)(-1.0,0.2)(-1.0,0)(-1.0,-0.2)(-0.9,-0.5)(-1.0,-1)
	}

\pscustom[linestyle=none, fillstyle=solid, fillcolor=gray!25]%
	{\pspolygon%
		(-1.5,-0.7)(-1.0,-0.6)(-0.5,-0.5)(-0.2,-0.6)(0,-0.6)(0.2,-0.6)(0.5,-0.5)(1.0,-0.6)(1.55,-0.7)%
		(1.4,-0.2)(1.4,0)(1.4,0.2)(1.5,0.5)%
		(1.0,0.6)(0.5,0.7)
	}

\pscustom[linestyle=none, fillstyle=solid, fillcolor=gray!35]%
	{\pspolygon%
		(-1.5,-0.7)(-1.0,-0.6)(-0.5,-0.5)(-0.2,-0.6)(0,-0.6)(0.2,-0.6)(0.5,-0.5)(1.0,-0.6)(1.55,-0.7)%
		(1.4,-0.2)(1.4,0)(1.4,0.2)(1.5,0.5)%
		(1.0,0.6)(0.5,0.7)(0.2,0.6)(0,0.6)(-0.2,0.6)(-0.5,0.7)(-1.0,0.6)(-1.5,0.5)%
		(-1.4,0.2)(-1.4,0)(-1.4,-0.2)(-1.3,-0.5)(-1.3,-0.7)
	}

		\rput(0,0){$\bullet$}

		\pscurve[linecolor=black, linewidth=1.5pt,arrows=<->, arrowscale=1.5](-2,0)(-1.5,-0.1)(-1,0)(-0.5,0.1)(-0.2,0)(0,0)(0.2,0)(0.5,0.1)(1,0)(1.5,-0.1)(2,0)

		\pscurve[linecolor=gray, linewidth=1.5pt,arrows=>-<, arrowscale= 1.5](0,-2) (-0.1,-1.5)(0,-1)(0.1,-0.5)(0,-0.2)(0,0)(0,0.2) (0.1,0.5)(0,1)(-0.1,1.5)(0,2)

	\pscurve[linecolor=gray](0.2,-1)(0.3,-0.5)(0.2,-0.2)(0.2,0)(0.2,0.2)(0.3,0.5)(0.2,1)
	\pscurve[linecolor=gray](0.4,-1)(0.5,-0.5)(0.4,-0.2)(0.4,0)(0.4,0.2)(0.5,0.5)(0.4,1)
	\pscurve[linecolor=gray](0.6,-1)(0.7,-0.5)(0.6,-0.2)(0.6,0)(0.6,0.2)(0.7,0.5)(0.6,1)
	\pscurve[linecolor=gray](0.8,-1)(0.9,-0.5)(0.8,-0.2)(0.8,0)(0.8,0.2)(0.9,0.5)(0.8,1)
	\pscurve[linecolor=gray](1.0,-1)(1.1,-0.5)(1.0,-0.2)(1.0,0)(1.0,0.2)(1.1,0.5)(1.0,1)
	\pscurve[linecolor=gray](1.25,-0.8)
(1.3,-0.5)(1.2,-0.2)(1.2,0)(1.2,0.2)(1.3,0.5)(1.2,0.75)
	\pscurve[linecolor=gray]
(1.5,-0.7)(1.4,-0.2)(1.4,0)(1.4,0.2)(1.5,0.5)

	\pscurve[linecolor=gray](-0.2,-1)(-0.1,-0.5)(-0.2,-0.2)(-0.2,0)(-0.2,0.2)(-0.1,0.5)(-0.2,1)
	\pscurve[linecolor=gray](-0.4,-1)(-0.3,-0.5)(-0.4,-0.2)(-0.4,0)(-0.4,0.2)(-0.3,0.5)(-0.4,1)
	\pscurve[linecolor=gray](-0.6,-1)(-0.5,-0.5)(-0.6,-0.2)(-0.6,0)(-0.6,0.2)(-0.5,0.5)(-0.6,1)
	\pscurve[linecolor=gray](-0.8,-1)(-0.7,-0.5)(-0.8,-0.2)(-0.8,0)(-0.8,0.2)(-0.7,0.5)(-0.8,1)
	\pscurve[linecolor=gray](-1.0,-1)(-0.9,-0.5)(-1.0,-0.2)(-1.0,0)(-1.0,0.2)(-0.9,0.5)(-1.0,1)
	\pscurve[linecolor=gray](-1.1,-0.8)(-1.1,-0.7)%
(-1.1,-0.5)(-1.2,-0.2)(-1.2,0)(-1.2,0.2)(-1.3,0.5)(-1.25,0.7)
	\pscurve[linecolor=gray](-1.3,-0.7)%
(-1.3,-0.5)(-1.4,-0.2)(-1.4,0)(-1.4,0.2)(-1.5,0.5)

	\pscurve[linecolor=black, linewidth=1.2pt](-1.5,0.1)(-1.0,0.2)(-0.5,0.3)(-0.2,0.2)(0,0.2)(0.2,0.2)(0.5,0.3)(1.0,0.2)(1.5,0.1)
	\pscurve[linecolor=black, linewidth=1.2pt](-1.5,0.3)(-1.0,0.4)(-0.5,0.5)(-0.2,0.4)(0,0.4)(0.2,0.4)(0.5,0.5)(1.0,0.4) (1.5,0.3)

	\pscurve[linecolor=black, linewidth=1.2pt] (-1.5,0.5)(-1.0,0.6)(-0.5,0.7)(-0.2,0.6)(0,0.6)(0.2,0.6)(0.5,0.7)(1.0,0.6) (1.5,0.5)

	\pscurve[linecolor=black](-1.25,0.7)(-1.0,0.8)(-0.5,0.9)(-0.2,0.8)(0,0.8)(0.2,0.8)(0.5,0.9)(1.0,0.8)(1.2,0.75)
	\pscurve[linecolor=black](-1.0,1.0)(-0.5,1.1)(-0.2,1.0)(0,1.0)(0.2,1.0)(0.5,1.1)(1.0,1.0)

	\pscurve[linecolor=black, linewidth=1.2pt](-1.5,-0.3)(-1.0,-0.2)(-0.5,-0.1)(-0.2,-0.2)(0,-0.2)(0.2,-0.2)(0.5,-0.1)(1.0,-0.2)(1.5,-0.3)
	\pscurve[linecolor=black, linewidth=1.2pt](-1.5,-0.5)(-1.0,-0.4)(-0.5,-0.3)(-0.2,-0.4)(0,-0.4)(0.2,-0.4)(0.5,-0.3)(1.0,-0.4)(1.5,-0.5)
	\pscurve[linecolor=black, linewidth=1.2pt](-1.5,-0.7)(-1.0,-0.6)(-0.5,-0.5)(-0.2,-0.6)(0,-0.6)(0.2,-0.6)(0.5,-0.5)(1.0,-0.6)(1.55,-0.7)
	\pscurve[linecolor=black](-1.1,-0.8)(-1.0,-0.8)(-0.5,-0.7)(-0.2,-0.8)(0,-0.8)(0.2,-0.8)(0.5,-0.7)(1.0,-0.8)(1.25,-0.8)
	\pscurve[linecolor=black](-1.0,-1.0)(-0.5,-0.9)(-0.2,-1.0)(0,-1.0)(0.2,-1.0)(0.5,-0.9)(1.0,-1.0)

	\pscurve[linestyle=dashed](-1.0,1.0)(-1.25,0.7)(-1.5,0.5)
	\pscurve[linestyle=dashed](1.0,1.0)(1.2,0.75)(1.5,0.5)
	\pscurve[linestyle=dashed](-1.3,-0.7)(-1.1,-0.8)(-1.0,-1)
	\pscurve[linestyle=dashed](1.55,-0.7)(1.25,-0.8)(1.0,-1)
	
	\uput[ul](-1,1){\qs{1}}
	\uput[ur](1,1){\qs{2}}
	\uput[dr](1,-1){\qs{3}}
	\uput[dl](-1,-1){\qs{4}}
	
	\uput[l](-1.5,0.5){\fof{\qs{1}}}
	\uput[r](1.5,0.5){\fof{\qs{2}}}
	\uput[r](1.5,-0.7){\fof{\qs{3}}}
	\uput[l](-1.5,-0.7){\fof{\qs{4}}}
	
	\rput(-1.5,1){\cGs{1}}\psline[linestyle=dotted](-1.5,1)(-1.25,0.7)
	\rput(1.5,1){\cGs{2}}\psline[linestyle=dotted](1.5,1)(1.25,0.75)
	\rput(1.5,-1){\cGs{3}}\psline[linestyle=dotted](1.5,-1)(1.25,-0.8)
	\rput(-1.5,-1){\cGs{4}}\psline[linestyle=dotted](-1.5,-1)(-1.1,-0.8)
	
	\rput(0.5,1.5){$\boldsymbol{\Eua}$}\psline[linestyle=dotted](0.5,1.5)(0.5,1.1)
	\rput(-0.5,-1.5){$\boldsymbol{\Eub}$}\psline[linestyle=dotted](-0.5,-1.5)(-0.5,-1.0)
	
	\rput(2.1,0.2){$\boldsymbol{\fof{\Eua}}$}\psline[linestyle=dotted](1.9,0.2)(1.4,0.2)
	\rput(-2.1,0.2){$\boldsymbol{\fof{\Eub}}$}\psline[linestyle=dotted](-1.9,0.2)(-1.4,0.2)

	\end{pspicture}\caption{\Vso}
\label{fig:Vso}
\end{center}
\end{figure}

The corresponding region \tVso{} for \A{} is defined by the inequalities
\begin{align*}
	\abs{xy}&\leq 1\\
	x^{2}&\leq 1\\
	y^{2}&\leq 1.
\end{align*}

We note that the two components $V^{u}_{i}$, $i=1,2$ of $\Vso\setminus\interior{\fof{\cR}}$ are \cFut-boxes, 
those of $\Vso\setminus\interior{\cR}$ ($V^{s}_{i}$) are \cFst-boxes, and all four 
are fundamental domains for $f$.  

\begin{lemma}\label{lem:honV}
	Let $V=\bigcup_{k\in\Integers}\toof{f}{k}{\Vso}$ and $\tV=\bigcup_{k\in\Integers}\toof{\A}{k}{\tVso}$.
	Then there is a homeomorphism $h$ conjugating $f|V$ with $\A|\tV$.
\end{lemma}
\begin{subproof}{\refer{lem}{honV}}
	First, we define $h$ on $X$: for $i=1,...,4$, 
	the interval \clint{\ps{i}}{\fof{\ps{i}}} is the one-dimensional analogue of a fundamental
	domain for $f$ (it is an interval abutting its $f$-image) and the analogue of \refer{rmk}{fundom}
	allows us to define a conjugacy between each separatrix of $f$ and the corresponding separatrix of \A.
	Since the orbit of each \ps{i} converges monotonically to the fixedpoint in one of the time directions,
	this definition, together with taking the fixedpoint to the origin, defines a homeomorphism $h$ taking $X$
	to the union of the two axes in \Realstwo, conjugating $f$ with \A{} there.
	
	Next, we use the foliations \cFst{} and \cFut{} to define a coordinate system on $V$: every point 
	$\vx\in\Vso$ is the (unique) point of intersection of the \cFst-arc through a point \xof{\vx} on the horizontal 
	arc in $X\cap\Vso$  with the \cFut-arc through some point \yof{\vx} on the vertical arc in $X\cap\Vso$;
	then the action of $f$ extends this property to all of $V$. 
	We define \hof{\vx}  to be the point $(\hof{\xof{\vx}},\hof{\yof{\vx}})$.	
	Note that  
	the images of any transversal to one of the separatrices 
	have as their limit set both of the "other" separatrices, together with the fixed point;
	this guarantees that the separate definitions of $h$ generated by the fundamental
	domains $V^{u}_{i}$ and $V^{s}_{i}$ agree on overlaps and have the right limit behavior at $X$. 
\end{subproof}

Finally, we extend the definition of $h$ on $V$ to each of the four components of the complement of $V$.
Each such component is a component of the complement of one of the leaves \cGs{\qsi} in the quadrant \Qsi.
By \refer{rmk}{quad}, there is a homeomorphism \vphi{} of \Qsi{} with the upper half plane that takes \cG-lines to
horizontal lines. The restriction of this homeomorphism to the union \pQsi{} of \cGs{\qsi} with the component of the
complement of $V$ which it bounds maps onto a closed half-plane, with \cGs{\qsi} going to the bounding horizontal line.
Let \Ts{0} be the \vphi-preimage of the vertical ray through some point \vphiof{p} on the boundary of the half plane: this is a global
cross section to the foliation \cG{} in \pQsi, as is its $f$-image $\Ts{1}=\fof{\Ts{0}}$. Let $S$ be the strip in \pQsi{} bounded by
the two cross-sections \Ts{0} and \Ts{1} together with the \cG-arc \Gclint{p}{\fof{p}}. 

The \vphi-image of \Ts{1} is not \emph{a priori} a vertical line; however, there is a homotopy of the plane, moving images of 
points along horizontal lines, which fixes $\vphiof{\Ts{0}\cup\Gclint{p}{\fof{p}}}$ and moves \vphiof{\Ts{1}} to a vertical line.
Composing $\vphi|S$ with this homotopy, we have a homeomorphism taking $S$ (which is a fundamental domain for $f$)
to a fundamental domain for the horizontal translation in the half plane. Applying \refer{rmk}{fundom}, we can extend this
homeomorphism to a conjugacy $h$ between $f|\pQsi$ and the horizontal translation in the half plane, which agrees with the 
previous definition of $h$ on \cGs{i}.
There is an easy corresponding conjugation of \A{} restricted to one of the components $\tQsi$ of the complement of \tV{}
and the horizontal translation on a half plane.  Composing the inverse of this conjugation with the one above gives
a homeomorphism between $f|\pQsi$ and $\A|\tQsi$ which agrees with the conjugacy $h|V$, defined previously, on the
common boundary.  

This proves 
\begin{prop}\label{prop:Mendes2}
	If $f$ is the time-one map of a \Cr{1} flow on \Realstwo{} with a single fixedpoint and has an Anosov structure, then it is topologically
	conjugate to the linear hyperbolic automorphism \A.
\end{prop}

In light of \refer{thm}{MendesThm}, \refer{cor}{Mendes1} and \refer{prop}{Mendes2} together prove Theorem A.

\bibliography{noncompact}
\bibliographystyle{amsplain}

\end{document}